\documentclass[12pt]{article} 
\title{\bf Computational Tools for the Shot Noise with Random Amplitude }
\author{Jean-Fran\c{c}ois Chamayou\thanks{Laboratoire de
Statistique et Probabilit\'es, Universit\'e Paul Sabatier, 31062 Toulouse, France. e-mail: \texttt{jfchamay@math.univ-toulouse.fr}}}
\date{\today}
\usepackage[cp850]{inputenc}
\usepackage[T1]{fontenc}
\usepackage[francais]{babel}
\usepackage{graphicx}

\def\a{{\bf a}}
\def\b{{\bf b}}
\def\A{{\bf A}}
\def\B{{\bf B}}

\begin{document}
\maketitle
\section{Abstract/Introduction}
 The following random recurrency:\\
$$
W_{n+1} = U_{n+1} ( W_n + \Lambda_{n+1} )
$$
where $W_{-1}=0.$ is known to be associated with the shot noise (see scheme 1):
$$
W_{t} = \sum_{0<t_k<t}  \Lambda_{k} e^{-(t-t_k)}
$$
where the $t_k$ are the dates of a Poisson process, the $ U_i
 ;\  i = 0,1,... $ are independent uniform variables.
 This can be extended to the triggered shot noise ( cf. ref. [7] , see scheme 2):
$$
W_{t} = \sum_{0<t_{lk}<t}  \Lambda_{lk}  e^{-(t-t_{lk})} , l=1,2,3,...
$$ 
The random amplitudes $\Lambda_i,\  i = 0,1,... $ are positive Dufresne independent variables which Laplace transform of the density is:\\
$$\,  _{p}F_{q} (\a;\b;-s) , p+1 \le q $$
or variables $ \lambda_c (\a;\b) $ independent which characteristic function is (cf. ref. [9]):\\
$$\,  _{2p+1}F_{2p} (c,\a,\a+\frac{1}{2};\b,\b+\frac{1}{2};-s^2) $$
or more generally even density variables which characteristic function (i.e. the Fourier transform) is:\\
$$\,  _{p}F_{q} (\a;\b,;-s^2) , p+1 \le q $$\\
the generalized hypergeometric function, with the sequences of parameters $p,q$; $ \a = a_1,...,a_p > 0 ; \b = b_1,...,b_q > 0 $
 defined for $-1< s< 1$:
\begin{eqnarray*}
\,   _{p}F_{q} (\a;\b;-s^2)=\sum_{n=0}^{\infty}\frac{ (\a)_n }{(\b)_n}\frac{(-s^{2})^{n}}{n!}.
\end{eqnarray*}
where $(a)_s=\frac{\Gamma(a+s)}{\Gamma(a)}$ is the Pochhammer symbol and $ \Gamma$ the Gamma function.\\
The compact notation used means for the products $ (\a)_n = (a_1)_n ... (a_p)_n $ and $ (-) $ an empty set of parameters.
Our aim is to provide tools to compute the density of the stationnary behaviour of the random recurrency i.e. of the classical or triggered shot noise. Special cases have already been studied

This work take advantage of previous results for different type of laws with the following transforms: \\
  Laplace\\
\begin{enumerate}
\item deterministic: $ \, _0F_0(-;-;-s)= \exp{(-s)} $, ref. [8],\\
\item Gamma: $ \, _1F_0(a;-;-s)= \frac{1}{(1+s)^a}, a>0$, ref. [7],\\
\item Beta: $ \, _1F_1(a;b;-s), b>a>0$,\\
\end{enumerate}
or Fourier:\\
\begin{enumerate}
\item Laplace: $ \, _1F_0(a;-;-s^2)= \frac{1}{(1+s^2)^a} , a>0$, ref. [7],\\
\item Cauchy: $ \, _0F_0(-;-;-|s|)= \exp{(-|s|)} $, ref. [7],\\
\item Standard Normal: $ \, _0F_0(-;-;-s^2/2)= \exp{(-s^2/2)} $, ref. [7],\\
\item Arcsine:  $ \, _0F_1(-;1;-s^2/4)= J_0(s) $, ref. [6,7], \\
\item Bernoulli: $ \, _0F_1(-;\frac{1}{2};-s^2/4)= \cos{(s)} $, ref. [5],[7], \\
\end{enumerate}
We intend to illustrate the general case with the particular example combinating the 2 last ones: sum of Arcsine and Bernoulli variables, i.e.: 
$$\, _0F_1(-;\frac{1}{2};-s^2/4) \, _0F_1(-;1;-s^2/4)= \cos{(s)} J_0(s)= \, _2F_3(\frac{1}{4},\frac{3}{4};\frac{1}{2},\frac{1}{2},1;-s^2) $$
see [20] 7.15.1.3 p 609.

\section{Shot Noise}
The stationnary density is obtained from a first order differential equation giving the transform (Laplace, Fourier) taking into account the characteristic function of the random amplitude $ g(s) $ (see [8] )
\begin{eqnarray*}
h(s) = \frac{C}{s} \exp(-\int_{s}^{\infty} \frac{g(\xi)}{\xi} d\xi)
\end{eqnarray*}
where $ C $ is the normalization factor.
The density of probability is given after the formal inversion of the Fourier transform for our general case by:
\begin{eqnarray*}
f(x) = \frac{C}{\pi} \int_{0}^{\infty} \cos(xt) \exp(-\int_{t}^{\infty} \frac{\,  _{p}F_{q} (\a;\b;-\xi^2)}{\xi} d\xi)
\frac{dt}{t} 
\end{eqnarray*}
the normalization factor is given by (see the appendix for the proof):
\begin{eqnarray*}
C = \exp(\frac{1}{2}(-\gamma-\psi(\a)+\psi(\b)))
\end{eqnarray*}
where $ \psi(a) $ is the logarithmic derivative of the $ \Gamma(a) $ function, $ \gamma $  the Euler constant, the compact notation used means for the sums: $ \psi(\a) = \psi(a_1) +...+\psi(a_p) $.

The previous papers required the use of known special functions: exponential, sine/cosine (hyperbolic),Bessel(modified)Integral, we introduce now the generalized hypergeometric Integral:
\begin{eqnarray*}
\, _{p}Ti_{q} (\a;\b;-x^2)  =2 \int_{x}^{\infty} \frac{ \,  _{p}F_{q} (\a;\b;-\xi^2)}{\xi} d\xi \\
=-\gamma-\ln(x)-\psi(\a)+\psi(\b) \\
+\,_{p}Tin_{q}(\a;\b;-x^2) 
\end{eqnarray*}
(see the appendix for the computation of the integration constant)
where the complementary function is:
\begin{eqnarray*}
\, _{p}Tin_{q} (\a;\b;-x^2)  =2 \int_{0}^{x} \frac{(1-\,  _{p}F_{q} (\a;\b;-\xi^2))}{\xi} d\xi
= \\
 - \sum_{k=1}^{\infty}  \frac{(\a)_{k}}{(\b)_{k} k k!}
 (-x^{2})^{k}.
\end{eqnarray*}

The asymptotical behavior for $ x \rightarrow \infty $ is given traditionally by integration by parts (for not integer $ \b $):\\

\begin{eqnarray*}
 2 \int_{x}^{\infty} \frac{ \,  _{p}F_{q} (\a;\b;-\xi^2)}{\xi} d\xi = \\
\sum_{k=1}^{\infty} \frac{(\b)_{-k} (k-1)! (-1)^{(k-1)}}{ (\a)_{-k} x^{2 k}}\,  _{p}F_{q} (\a-k;\b-k;-x^2).
\end{eqnarray*}

We could have introduced even more general function like Meijer function 
where the hypergeometric function is replaced by its extension (see [20] III,8.2.2.3 page 618):
\begin{eqnarray*}
G^{m,n}_{p,q} \left[ \xi \left| \begin{array} {c} \a \\ \b \end{array} \right| \right] = \\
\sum_{k=1}^{m} \frac{ \prod_{j=1}^{* \, m} \Gamma(b_j-b_k) \prod_{j=1}^{n}  \Gamma(1-a_j+b_k)}
{\prod_{j=m+1}^{q} \Gamma(1-b_j+b_k) \prod_{j=n+1}^{p} \Gamma(a_j-b_k)}\\
\,_pF_{q-1} ( 1-a_1+b_k,...,1-a_p+b_k ;
 1-b_1+b_k,.*.,1-b_q+b_k;(-1)^{p-m-n} \xi ) 
\end{eqnarray*}
where the normalization gives:
\begin{eqnarray*}
G^{m,n}_{p,q} \left[ 0 \left| \begin{array} {c} \a \\ \b \end{array} \right| \right] = 
\sum_{k=1}^{m} \frac{ \prod_{j=1}^{* \, m} \Gamma(b_j-b_k) \prod_{j=1}^{n}  \Gamma(1-a_j+b_k)}
{\prod_{j=m+1}^{q} \Gamma(1-b_j+b_k) \prod_{j=n+1}^{p} \Gamma(a_j-b_k)}\\
\end{eqnarray*}

\subsection{Asymptotical Study}

An asymptotical study by the saddle point method of the stationary density is undertaken by considering the inverse Laplace transform:
$$
f(x) = \frac{1}{2i\pi} \int_{A_0} \exp(- sx +\int_{0}^{s}
\frac{\,    _{p}F_{q} (\a;\b;\xi^2)-1}{\xi} d\xi) ds
$$
where $A_0 $ is a vertical path in the $s$ plane through the saddle point $s_0$ if
$$ \Phi(s) = - sx +\int_{0}^{s} \frac{\,  _{p}F_{q} (\a;\b;\xi^2) -1 }{\xi}d\xi,$$ 
the saddle points 
 $s_0$ are the roots of
$$
 \Phi^{'}(s) = -x+ \frac{\,  _{p}F_{q} (\a;\b;s^2)-1}{s}.
$$
 Then the density $f(x)$ can be approximated by:
$$
 f(x)= \frac{\exp(\Phi(s_0))}{\sqrt{2\pi \Phi^{"}(s_0)}}, x \rightarrow \infty
$$
For $x$ sufficiently large $\Phi^{'}(s) $ has only one root 
$s_0(x)$ tending to infinity with $x. $ This fact can be seen using the asymptotical representation of $ \, _{p}F_{q} (\a;\b;s^2) $
The calculation of the roots of $$ 1+x s_0 =  \,  _{p}F_{q} (\a;\b;s_0^2) $$
 gives an idea of the asymptotic behavior, first with:  $ s_0^{(1)} = O(\ln (x))$ 
 
\subsection{Our Favorite Example}

We handle the following random recurrency:
\begin{eqnarray*}
W_{n+1} = U_{n+1} ( W_n + cos(\pi V_{n+1}) + \Delta_{n+1} )
\end{eqnarray*}
where:
\begin{eqnarray*}
W_0 = U_0 (  cos(\pi V_0) + \Delta_0) ;  n = 0,1,...
\end{eqnarray*}
 the density of $W_0$ is the following:
$ f_0(x)= \frac{1}{2 \pi}\sqrt{\frac{2-|x|}{|x|}} {\bf 1}_{[-2,2]}(x)$\\
$f_1(x) $ the density of $ W_1 $ is given from the intermediary step: calculation of 
$ g_1(y) $ the density of $ Y = W_0 + (  cos(\pi V_1) + \Delta_1) $ obtained using complete
elliptic functions $ { \bf K } $ :
\begin{eqnarray*}
g_1(y) = 
\frac{1}{\pi^2}(2 \, {\bf K}(1-\frac{y^2}{4}) {\bf 1}_{[-2,2]}
+(2-\frac{y}{2})\, {\bf K }((2-\frac{y}{2})\frac{y}{2}) {\bf 1}_{[-4,4]})(y) 
 \end{eqnarray*}
 Then $ f_1 $ can be expressed through the integral:
\begin{eqnarray*}
\int_0^{\sqrt{1-(x-1)^2}}\, \frac{{\bf K}(\xi)}{\xi}\,
 (\sqrt[5]{1-\xi^2}-\frac{1}{\sqrt[5]{1-\xi^2}})^2 d\xi 
 \end{eqnarray*}
 integrated by series expansion see Prudnikov , vol III page 36. The comparisons between the densities and the Monte-Carlo simulations are presented in figure 1 for $f$ and figure 2 for $g$ (only on the positive abscissae due to parity).
 
The previous recurrency is associated to the shot noise: 
\begin{eqnarray*}
W_{t} = \sum_{0<t_k<t}  ( cos(\pi V_{n+1}) + \Delta_{n+1} ) e^{-(t-t_k)}
\end{eqnarray*}
where the $ t_k$ are the dates of a Poisson process, the 
$ U_i , V_i ; i = 0,1,... $ are independent uniform random variables and the $ \Delta_i , i = 0,1,... $ 
independent Bernoulli variables $ \frac{1}{2} ( \delta_{+1} + \delta_{-1} ) $.\\
This random amplitude can  be considered as the product of the same Bernoulli
by a Beta $ \beta_{\frac{1}{2},\frac{1}{2}}^{(1)} $ multiplied by 2, this fact is also shown by the Carlson identity [4] entry (11)$p=1, a=\frac{1}{2}, q=1, b=1$:

\begin{eqnarray*}
\,_{2 p}F_{2q+1}( \frac{\a}{2},\frac{\a+1}{2};\frac{1}{2},\frac{\b}{2},\frac{\b+1}{2};-s2)=\\
\frac{1}{2}(\,_{p}F_{q}( \a;\b;i 2^{1+q-p}s)+\,_{p}F_{q}( \a;\b;-i 2^{1+q-p}s))
\end{eqnarray*}

A numerical computation of the stationary density is tractable by splitting the integration domain of the integral and negleting the minor integration part (see appendix for the normalization constant):\\
\begin{eqnarray*}
f(x) = \frac{2e^{-\gamma}}{\pi} \int_{0}^{\infty} cos(xt) exp(- \int_{t}^{\infty} \frac{cos(\xi) J_0 (\xi)}{\xi} d\xi) \frac{dt}{t} =
\end{eqnarray*}
\begin{eqnarray*}
\frac{2e^{-\gamma}}{\pi} [\int_{0}^{x_1} cos(xt) exp(- \int_{t}^{\infty} \frac{cos(\xi) J_0 (\xi)}{\xi} d\xi) \frac{dt}{t} \\
-\frac{ sin(xx_1)}{xx_1} exp(- \int_{x_1}^{\infty} \frac{cos(\xi) J_0 (\xi)}{\xi} d\xi) \\
+\int_{x_1}^{x_2} \frac{1-cos(t)J_0(t)}{xt} sin(xt) exp(- \int_{t}^{\infty} \frac{cos(\xi) J_0 (\xi)}{\xi} d\xi) \frac{dt}{t}]+O(\frac{1}{xx_2^2})
\end{eqnarray*}
where $O(\frac{1}{xx_2^2}) \sim O(h^2) , h $ is the integration step from the trapezoidal rule of this oscillating integral (see [11] for the numerical results).\\

 The computations are based on earlier results obtained with other types of laws (Arcsine or Bernoulli alone), they require now the calculations of the special function Cosine-Bessel Integral in accordance with the Agrest's papers [1,2,3] for the function Exponential-Bessel Integral.
We can choose between several expansions:
\begin{eqnarray*}
 cos(x) J_0 (x)
                & = &\sum_{k=0}^{\infty} \frac{(\frac{1}{2})_{2k}  }{( (2k)! )^2}(-4x^{2})^{k}.
\end{eqnarray*}
which gives the following moments: $ K_{2k} = \frac{(\frac{1}{2})_{2k} 2^{2k} }{(2k)!},K_{2k+1}= 0 ,  k=0,1,... $, or:
\begin{eqnarray*}
 cos(x) J_0 (x)
                & = &
\sum_{k=0}^{\infty} (-1)^k \frac{(4k)!}{ ((2k)!)^3} (\frac{x}{2})^{2k}
\end{eqnarray*}
found in [13] 5-24-30 p 92, and the hypergeometrical representation already established
\begin{eqnarray*}
 cos(x) J_0 (x) = \, _2F_3( \frac{1}{4},\frac{3}{4}; \frac{1}{2},\frac{1}{2},1 ; - x^2 ) = \sum_{n=0}^{\infty}\frac{(\frac{1}{4})_n (\frac{3}{4})_n }{(\frac{1}{2})_n (\frac{1}{2})_n}\frac{(-x^{2})^{n}}{(n!)^2}.
\end{eqnarray*}
which is more convenient for our purpose, we get:
\begin{eqnarray*}
CJi ( x )  = 
\int_{x}^{\infty}  \frac{cos(\xi) J_0 (\xi)}{\xi} d\xi =
\frac{1}{\pi} \int_{0}^{\pi} Ci(2 x \sin^{2}{(\frac{\theta}{2})}) d\theta =
 m_1  - Ln(x)+ CJin (x) 
\end{eqnarray*}
where $ Ci $ is the cosine integral (the integral uses the integral representation of the Bessel function)and the complementary function is: 
\begin{eqnarray*}
CJin (x) = \int_{0}^{x} \frac{(1 - cos(\xi) J_0 (\xi))}{\xi} d\xi =
-\sum_{n=1}^{\infty}\frac{(\frac{1}{4})_n (\frac{3}{4})_n }{(\frac{1}{2})_n (\frac{1}{2})_n}\frac{(-x^{2})^{n}}{2n (n!)^2}
\end{eqnarray*}
$ m_1 $ being the integration constant computed in the appendix,
the constant is given by:
$$ m_1=  \frac{1}{2} ( - \gamma - \psi(\frac{1}{4}) - \psi(\frac{3}{4}) + 2 \psi(\frac{1}{2}) +\psi(1) ) =
 - \gamma + Ln(2) $$
(cf. ref. [21]  for the $ \psi $ values or the  pseudo code of its computation).\\
The asymptotical behavior for $ x \rightarrow \infty $ is given by an integration by parts (due to the integer parameter of the hypergeometric function) using the Hankel asymptotic expansion for $ J_0 $ see Abramovitz [0] page 364:
\begin{eqnarray*}
 \int_{x}^{\infty} \frac{cos(\xi) J_0 (\xi)}{\xi} d\xi \sim
 \frac{2-\frac{1}{12 x} +\frac{\sqrt{2} \cos{(2 x -\frac{\pi}{4})}}{4 x}}{\sqrt{\pi x}}
 +\circ(\frac{1}{x^{\frac{5}{2}}})
\end{eqnarray*}

\subsection{Asymptotical Study of the example}

The asymptotic behavior of the density is given from the saddle point method by considering the bilateral inverse Laplace transform:\\
\begin{eqnarray*}
f(x) = \frac{1}{2i\pi} \int_{A_0} exp(- sx +\int_{0}^{s} \frac{Ch(\xi) I_0 (\xi)-1}{\xi} d\xi) ds 
\end{eqnarray*}
where $A_0 $ is a vertical path in the $ s$ plane through the saddle point  $s_0$ \\
if $ \Phi(s) = - sx +\int_{0}^{s} \frac{Ch(\xi) I_0 (\xi)-1}{\xi} d\xi $ the saddle points  $s_0$ are the roots from:\\
\begin{eqnarray*}
 \Phi^{'}(s) = -x+ \frac{Ch(s) I_0 (s)-1}{s} 
\end{eqnarray*}
 then the density $ f(x) $ can be approximated by:\\
\begin{eqnarray*}
 f(x)= \frac{e^{\Phi(s_0)}}{\sqrt{2\pi \Phi^{"}(s_0)}}, x \rightarrow \infty 
\end{eqnarray*}
for $x$ sufficiently large $\Phi^{'}(s) $ has only one root $ s_0(x) $ tending to infinity with $x$,and using the following asymptotic representations for the modified Bessel function and the hyperbolic functions:\\
\begin{eqnarray*}
 I_n(s)\sim \frac{e^s}{\sqrt{2 \pi s}} , Ch(s) \sim Sh(s) \sim \frac{e^s}{2} 
\end{eqnarray*}
the calculation of the roots from $ 1+xs_0 = \frac{e^{2s_0}}{2\sqrt{2 \pi s_0}} $ gives an idea of the asymptotic behavior, first  $ s_0^{(1)} = O(\frac{1}{2} Ln(x))$ and $ s_0^{(2)} = \frac{1}{2} Ln(x) + \frac{3}{4} Ln(Ln(x)) + \frac{1}{4} Ln(\pi) , x \rightarrow \infty $\\
The analogous relations exist for the special function hyperbolic cosine modified Bessel integral for the asymptotic study with the help of the saddle point method:
\begin{eqnarray*}
CIi ( x )  = \gamma + Ln(x/2)  + CIin (x) 
\end{eqnarray*}
where:
\begin{eqnarray*}
CIin (x) = \int_{0}^{x} \frac{(Ch(\xi) I_0 (\xi)-1)}{\xi} d\xi =
\sum_{n=1}^{\infty}\frac{(\frac{1}{4})_n (\frac{3}{4})_n }{(\frac{1}{2})_n (\frac{1}{2})_n}\frac{(x^{2})^{n}}{2n (n!)^2}
\end{eqnarray*}
again the asymptotic behavior for $ x \rightarrow \infty $ is given by an integration   using the asymptotic expansion of $ I_0 $ see Abramovitz [0] page 377:
\begin{eqnarray*}
 \int \frac{Ch(\xi) I_0 (\xi)}{\xi} d\xi \sim
  \frac{\sqrt{2} \exp{(2 x)}}{8 x \sqrt{\pi x}}
 +\circ(\frac{\exp{(2 x)}}{x^{\frac{5}{2}}})
\end{eqnarray*}

\section{Triggered Shot Noise}

 \subsection{General Case}
 
 The approximated stationary density of the random recurrency is obtained from a  lth order differential equation giving the transforms (Laplace, Fourier) $ h(s) $ taking into account the characteristic function of the random amplitude $ g(s) $ (see [7] ), for l=2:
\begin{eqnarray*}
s^2 h''(s) + 3 s h'(s) + ( 1 - g(s) ) h(s) = 0 
\end{eqnarray*}
the $h$ expansion in the neighbourhood of zero is:
$ h(s) = \sum_{i=0}^{\infty} c_i s^i $, the coefficients are identified using the $ g(s) $
expansion, we get:
\begin{eqnarray*}c_0=1 , c_{2n+1}=0, 
c_{2n}=\frac{1}{4n(n+1)}\sum_{j=1}^{n-1}\frac{c_{2j} K_{2(n-j)}}{(2(n-j))!}
, n=1,2,...
\end{eqnarray*}
where the $ K_n $ are the moments of the random amplitude.
The iterative computation according to cf. [6] of $ h(s) $:
$$ s h_2(s) = C_1 (( Ln(s) + \gamma)(1+Ti^{(2)}(s))+2 Ti^{(3)}(s))+(C_2-C_1)(1+Ti^{(2)}(s))$$
is used for the calculation of the 2 constants $ C_1, C_2 $ from the comparison with $ \sum_{j=0}^{\infty} c_j s^j $ which allows the knowledge of the asymptotical behavior for the triggered case.\\
For the triggered shot noise we need to introduce the special function iterated generalized hypergeometric Integral function:
$$ \, _{p}Ti_{q}^{(n)} (\a;\b;-x^2)  = 2 \int_{x}^{\infty} \frac{ \,_{p}Ti_{q}^{(n-1)} (\a;\b;\xi^2)}{\xi} d\xi $$
Where:
$$\, _{p}Ti_{q}^{(0)} (\a;\b;-x^2) =2 \,  _{p}F_{q} (\a;\b;-x^2)$$
and
$$ \, _{p}Tin_{q}^{(n)} (\a;\b;-x^2) =
 -\sum_{k=1}^{\infty}   \frac{(\a)_{k} }{(\b)_{k} (k)^{n} k!}
 (-x^{2})^{k}.
$$\\
The semiconvergent asymptotic forms $ x \rightarrow \infty $ of the iterated exponential integral can be found in Kunstner [17], Van De Hulst [22] page 11 and of the iterated cosine integral in Hallen [12], for the iterated Bessel integral we get the following result from Bessel Integral $ Ji^{(1)}(x)$ approximation found in Petiau [19] Smith, page 254 :\\
\begin{eqnarray*}
Ji^{(2)}(x) = J_0(x) ( \frac{1}{x^2}-\frac{20}{x^4}+... +(-1)^{n-1} \frac{(2^n n!)( 2^{n-1} (n-1)!)(\sum_{k=1}^{n-1}\frac{1}{k} +\frac{1}{2n} )}{x^{2n}}+...)\\
+ J_1(x) ( \frac{4}{x^3}-\frac{96}{x^5} +...+(-1)^{n-1} \frac{(2^n n!)^2 (\sum_{k=1}^{n}\frac{1}{k} )}{x^{2n+1}}+...)
\end{eqnarray*}
\subsection{Our favorite Example}

According to the results in [7] on the triggered shot noise we get for the $ W_0 $ density:\\
$$ f_0(x)= \frac{1}{\pi}(\sqrt{\frac{2-|x|}{|x|}}- 
Arctg(\sqrt{\frac{2-|x|}{|x|}})){\bf 1}_{[-2,2]}(x)$$
figure 3 shows the comparison with a Monte-Carlo simulation of the process.\\
The link between iterated Integral functions and their complementary is given, cf. [14] for $ \psi^{(n)} $ using $ m_2 , m_3 $ (see appendix), by:
\begin{eqnarray*}
 m_2 = \frac{1}{2^2 2!} (\psi'(\frac{1}{4}) + \psi'(\frac{3}{4}) - 2 \psi'(\frac{1}{2}) ) = \pi^2/4 \\
 m_3=\frac{1}{2^3 3!}(2\psi''(1)-(\psi''(\frac{1}{4})+\psi''(\frac{3}{4}))+2\psi''(\frac{1}{2}))
=\\
\frac{1}{16.3}(-4\zeta(3)+4.28\zeta(3)-2(8-1)\zeta(3))=\zeta(3)/3 
\end{eqnarray*}
then:
\begin{eqnarray*}
CJi^{(2)}(x) = \frac{1}{2!} (\gamma +Ln(x/2))^2 + \pi^2/4 
-\sum_{n=1}^{\infty}\frac{(\frac{1}{4})_n (\frac{3}{4})_n }{(\frac{1}{2})_n (\frac{1}{2})_n}\frac{(-x^{2})^{n}}
{(2n)^2 (n!)^2}\\
CJi^{(3)}(x) = -\frac{1}{3!} (\gamma +Ln(x/2))^3 +\frac{\pi^2}{4} (\gamma+Ln(x/2))+\zeta(3)/3
-\sum_{n=1}^{\infty}\frac{(\frac{1}{4})_n (\frac{3}{4})_n }{(\frac{1}{2})_n (\frac{1}{2})_n}\frac{(-x^{2})^{n}}
{(2n)^3 (n!)^2}
\end{eqnarray*}
where the zeta function: 
$$ \zeta(3) = 1.2020569032$$ 

Since these functions are necessary to the computation according to cf. [6] of:
$$ s h_2(s) = C_1 (( Ln(s) + \gamma)(1+CJi^{(2)}(s))+2 CJi^{(3)}(s))+(C_2-C_1)(1+CJi^{(2)}(s)) + \o{1/(\sqrt{s})^9}$$
for the calculation of the 2 constants $ C_1, C_2 $ from the comparison with $ \sum_{j=0}^{\infty} c_j s^j $ which allows the knowledge of the asymptotical behavior for the triggered case, figure 4 shows the comparison (only for positive abscissae due to the parity).\\
The asymptotical behavior for $ x \rightarrow \infty $ is given by a double integration by parts using the Hankel asymptotic expansion for $ J_0 $ see Abramovitz [0] page 364:
\begin{eqnarray*}
 \int_{x}^{\infty} \frac{1}{\xi}  \int_{\xi}^{\infty} \frac{cos(\eta) J_0 (\eta)}{\eta} d\eta d\xi =
 \frac{4-\frac{1}{6 x} }{\sqrt{\pi x}} + +\circ(\frac{1}{x^{\frac{5}{2}}})
 \end{eqnarray*}
 The knowledge of the asymptotic behaviour is necessary to compute by inverse Fourier transform the density with short support from the characteristic function with extended support as shown in figure (5)
\subsection{Waiting Time Paradox }
It is necessary to take into account the waiting time correction for the triggered shot noise, (see [23],[7]) when $ l =2 $ for instance:
$$ f_{waiting-time}(x)= \frac{1}{4\pi}(3 \sqrt{\frac{2-|x|}{|x|}}- 
 2 Arctg(\sqrt{\frac{2-|x|}{|x|}})){\bf 1}_{[-2,2]}(x)$$
 
 The density of the product of the random amplitude by 3 independent uniforms variables  can also be calculated in closed form:
\begin{eqnarray*}
f_{0}(x)= \frac{2}{\pi}(\sqrt{\frac{2-|x|}{|x|}}- 
 Arctg(\sqrt{\frac{2-|x|}{|x|}})\\
 -\frac{1}{2} (\omega + \theta ) Ln(r) +\frac{1}{4} ( Cl_2 (2(\omega+\theta)) 
-(Cl_2 (2\omega)+Cl_2(2\theta)))){\bf 1}_{[-2,2]}(x)
\end{eqnarray*}
 where: $ r= \frac{1}{\sqrt{2|x|}} , \theta= Arctg(\sqrt{\frac{2-|x|}{|x|}})$ and
 $ \omega = Arctg(\frac{r sin(\theta)}{(1- r cos(\theta))})$\\
with the help of the Clausen Integral, see [18]:
 $$ Cl_2(\alpha) = -\sum_{k=1}^{\infty} \frac{sin( k\alpha)}{k^2} = 
 - \int_{0}^{\alpha} Ln(2 sin(\frac{\theta}{2})) d\theta $$
 Figure 5 shows the comparison of $ f_0 $ with the Monte-Carlo simulation of the random recurrency, the first term is sufficient, except for an accurate computation of the tail of the stationnary density.

 All the other numerical resuls and graphs can be found in the report [11].\\

\section{Appendix: Auxiliary Result}

The method to compute the integration constants of the iterated Integral functions follows the line of the proof in [16] pages 258-259.\\
{\bf Lemma 1:}\\

Let $ g(s) $ be the transform (Laplace or Fourier) of the random amplitude density, 
with moments : $ K_k, k=1,2,...$, we get:
$$ g^{(n)}(s) = \int_s^{\infty} g^{(n-1)}(\xi) \frac{d \xi}{\xi},
g^{(0)}(s) = g(s) $$
we get:
\begin{eqnarray*}
 g^{(n)}(s) =
M_{n-1}+(-1)^n \frac{ Ln^{n} (s)}{n!}
+(-1)^n \sum_{k=1}^{\infty} (-1)^{k-1} \frac{K_k s^k}{k^n k!}
+\sum_{l=0}^{n-2} \frac{(-1)^{(n-l-1)} Ln^{n-l-1}(s) M_l }{l! (n-l-1)! }
\end{eqnarray*}
where:
\begin{eqnarray*}
M_n=
\int_{0}^{\infty} g'(\xi) \frac{Ln^{n} (\xi)}{ n! }d \xi
\end{eqnarray*}
{\bf Proof:}\\
Integrating by parts:
$$ g^{(n)}(s) = 
\int_1^{\infty} g(s \xi) \frac{Ln^{n-1}(\xi)}{ (n-1)! }\frac{ d \xi}{ \xi} 
, n=1,2,...$$
Splitting the domain of integration we get for $g^{(n)}(s)$ :
\begin{eqnarray*}
\int_{\frac{1}{s}}^{\infty} g(s \xi) \frac{Ln^{n-1}(\xi)}{(n-1)!}\frac{d \xi}{\xi}
- \int_{0}^{\frac{1}{s}}(1- g(s \xi)) \frac{Ln^{n-1} (\xi)}{(n-1)!}\frac{d \xi}{\xi}\\
+\int_{1}^{\frac{1}{s}} \frac{Ln^{n-1} (\xi)}{ (n-1)! }\frac{ d \xi}{ \xi}
+\int_{0}^{1} (1-g(s \xi)) \frac{Ln^{n-1}(\xi)}{(n-1)!}\frac{ d \xi}{\xi}= \\
\int_{1}^{\infty} g(\xi) \frac{Ln^{n-1} (\xi/s)}{ (n-1)! }\frac{ d \xi}{\xi} 
-\int_{0}^{1}(1- g(\xi)) \frac{Ln^{n-1} (\xi/s)}{ (n-1)! }\frac{ d \xi}{\xi}\\
+(-1)^n\frac{ Ln^{n}(s)}{n!} + \sum_{k=1}^{\infty}(-1)^{k-1}\frac{K_k s^k}{k!}
\int_{0}^{1} x^{k-1} \frac{Ln^{n-1} (\xi)}{ (n-1)! }d \xi=\\
 \int_{1}^{\infty} g(\xi) \frac{Ln^{n-1} (\xi)}{ (n-1)! }\frac{ d \xi}{ \xi} 
-\int_{0}^{1}(1- g(\xi)) \frac{Ln^{n-1} (\xi) }{ (n-1)! }\frac{d \xi}{\xi}\\
+(-1)^n\frac{Ln^{n}(s)}{n!}
+(-1)^n\sum_{k=1}^{\infty}(-1)^{k-1}\frac{K_k s^k}{k^n k!}\\
+\sum_{l=0}^{n-2} \frac{(-1)^{(n-l-1)} Ln^{n-l-1}(s)}{l! (n-l-1)! }
[ \int_{1}^{\infty} g(\xi) \frac{Ln^{l} (\xi) d \xi}{\xi} 
-\int_{0}^{1}(1- g(\xi)) \frac{Ln^{l} (\xi) d \xi}{\xi}]=\\
M_{n-1}
+(-1)^n \frac{ Ln^{n} (s)}{n!}
+(-1)^n \sum_{k=1}^{\infty} (-1)^{k-1} \frac{K_k s^k}{k^n k!}
+\sum_{l=0}^{n-2} \frac{(-1)^{(n-l-1)} Ln^{n-l-1}(s) M_l }{l! (n-l-1)! }
\end{eqnarray*}
The integrals $ M_n$ are also obtained by integration by parts:
\begin{eqnarray*}
M_n = \int_{1}^{\infty} g(\xi) \frac{Ln^{n} (\xi)}{ n! }\frac{ d \xi}{\xi} 
-\int_{0}^{1}(1- g(\xi)) \frac{Ln^{n} (\xi) }{ n! }\frac{d \xi}{\xi}=
\int_{0}^{\infty} g'(\xi) \frac{Ln^{n} (\xi)}{ n! }d \xi
\end{eqnarray*}\\
$\hfill  \ddag  $\\
Generally this integral is unknown in closed form for any $ g(s) $, more likely we will obtain the following one:
\begin{eqnarray*}
M_n=\frac{1}{n!}[\frac{d^n}{ds^n}(K(s))]_{s=0}=
\frac{1}{n!}[\frac{d^n}{ds^n}(\int_{0}^{\infty}g'(\xi)\xi^s d\xi)]_{s=0}
\end{eqnarray*}

This is true for instance for the Fox $ H $ function, a very general special function
(see [20] III paragraph 8.3 pp626):
\begin{eqnarray*}
g(s) =  H^{m,n}_{p,q} \left[ s \left| \begin{array}{c}  (\a, \A) \\(\b, \B) \end{array}  \right| \right]
\end{eqnarray*}
where(cf. [20]III,8.3.2.15 page 629)for the differentiation of the Fox function ):
\begin{eqnarray*}
s g'(s) =  H^{m,n+1}_{p+1,q+1} \left[ s \left| \begin{array}{c} (0,1);(\a, \A) \\ (\b, \B);(1,1)  \end{array}\right| \right]      
\end{eqnarray*}

From the definition, using Mellin-Barnes integrals (see [20]III,8.3.1.1 page 626) we get for the integral $ K $:
\begin{eqnarray*}
K(s) =  \int_{0}^{\infty} \xi^{s-1} H^{m,n+1}_{p+1,q+1} 
\left[ \xi \left| \begin{array}{c}  (0,1); (\a, \A) \\  (\b, \B);(1,1) \end{array} \right| \right] d\xi= \\
-s \frac{\prod_{j=1}^m \Gamma(b_j-B_js) \prod_{j=1}^n \Gamma(1-a_j-A_js)}
			 {\prod_{j=m+1}^q \Gamma(1-b_j-B_js) \prod_{j=n+1}^p \Gamma(a_j-A_js)}
\end{eqnarray*}

Returning to our general case:\\
{\bf Lemma 2:}\\
For
$$ g(s) = \, _pF_q(\a;\b;-s^2)$$
the integrals:
\begin{eqnarray*}
M_n= \frac{1}{n!}[\frac{d^n}{ds^n}(\int_{0}^{\infty}g'(\xi)\xi^s d\xi)]_{s=0} , n=1,2...
\end{eqnarray*}
are given by the recurrency:
$$M_n = m_n + m_1 M_{n-1} , M_0 = 0 , n=1,2,...$$
where:
$$m_k =\frac{1}{2^k k!} (\psi^{(k-1)}(1)+(-1)^k(\psi^{(k-1)}(\a)-\psi^{(k-1)}(\b))), k=1,2,...$$ 
 the $ \psi^{(k)}$ are the derivative of the $ \psi $ function.\\
 {\bf Proof.}\\
( cf. [20] for the differentiation of the hypergeometric function ):
$$ g'(s) = -2s \frac{(\a)_1 }{(\b)_1} \,_pF_q(\a+1;\b+1;-s^2)$$
then from the definition of the generalized hypergeometric functions using Mellin-Barnes integrals 
cf. [20] 7.3.4.12 page 438. we get:
$$ K(s) = \frac{(\a)_1 }{(\b)_1}\int_0^{\infty} \, _pF_q(\a+1;\b+1;-\xi) 
\xi^{s/2} d\xi = \frac{(\a)_{-s/2} \Gamma(1+s/2)}{(\b)_{-s/2}}
$$
from which:
\begin{eqnarray*} 
K'(s)= \frac{1}{2} (\frac{(\a)_{-s/2} \Gamma'(1+s/2)}{(\b)_{-s/2}}
+\frac{(\a)_{-s/2} \Gamma(1+s/2)\Gamma(\b) \Gamma'(\b-s/2)}{\Gamma^2(\b-s/2)} 
-\frac{\Gamma(1+s/2)\Gamma'(\a-s/2)}{\Gamma(\a)(\b)_{-s/2}})
\end{eqnarray*}
then:
\begin{eqnarray*}
M_1= m_1 = K'(0) =  \frac{1}{2}( \psi(1)+\psi(\b)-\psi(\a))
\end{eqnarray*}
 However to continue it will be more convenient to introduce the expansion of $ Ln(K(s)) $ using the expansions of the logarithm of the $ \Gamma $ function, cf. [13] 54-111 p 358 where the Riemann Hurwitz function is given in terms of derivative of the $ \psi $ function:
$$ \zeta(k,\a)=\frac{(-1)^k}{(k-1)!} \psi^{(k-1)}(\a)$$ 
we get:
\begin{eqnarray*}
Ln(K(s)) = Ln(\frac{\Gamma(\a-s/2)}{\Gamma(\a)})
-Ln(\frac{\Gamma(\b-s/2)}{\Gamma(\b)})
+Ln(\frac{\Gamma(1+s/2)}{\Gamma(1)})=\\
\frac{1}{2}(\psi(1)-\psi(\a)+\psi(\b)) t/2 +
\sum_{k=2}^{\infty}\frac{t^k}{2^k k!} (\psi^{(k-1)}(1)+(-1)^k(\psi^{(k-1)}(\a)-\psi^{(k-1)}(\b))) 
\end{eqnarray*}
finally:
$$m_k =\frac{1}{2^k k!} (\psi^{(k-1)}(1)+(-1)^k(\psi^{(k-1)}(\a)-\psi^{(k-1)}(\b))), k=1,2,...$$ 
the $ M_n $ are obtained from the  $ m_k$ using the application of the method presented in [10],[14] :
$$M_n = m_n + m_1 M_{n-1} , M_0 = 0 , n=1,2,...$$
$\hfill  \ddag  $\\
 
 \newpage
\section{References:}

\vspace{4mm}\noindent[0] \textsc{Abramovitz M., Stegun I.A.} 
{\it Handbook of Mathematical Functions }{\bf [1964]}, National Bureau of Standards.

\vspace{4mm}\noindent[1] \textsc{M.M. Agrest}
{\it "Generalization of some relations for the Bessel Integral Functions"}
{\bf Bull. Acad. Sci. Georgian SSR }126,2{\bf [1987]} 241-244.

\vspace{4mm}\noindent[2] \textsc{M.M. Agrest, T.S. Chachibaya}
{\it " Expansions of Bessel Exponential-Integral Functions and related Functions"}
{\bf Soviet Math. ( Iz VUZ )} 32,4{\bf [1988] }1-13.

\vspace{4mm}\noindent[3] \textsc{M.M. Agrest}
{\it " Certain Relations for Exponential-Integral Bessel Functions of genus 1 "}
{\bf Soviet Math. ( Iz VUZ)} 35,6 {\bf[1991]} 67-69. 

\vspace{4mm}\noindent[4] \textsc{B. C. Carlson}
{\it " Some extensions of Lardner's relations between $\, _0F_3$ and Bessel functions "}
{\bf SIAMJ. Math. Anal.} 1 {\bf[1970]} 232-242,MR41:3819.

\vspace{4mm}\noindent[5] \textsc{J-F. Chamayou}
{\it " Numerical evaluation of a solution of a special mixed type differential difference equation "}
{\bf Calcolo} 15 {\bf[1978]} 395-414. 

\vspace{4mm}\noindent[6] \textsc{J-F. Chamayou}
{\it " Mod\`{e}le de bruit de grenaille trigonom\'etrique "}
{\bf Journ\'ees de Statistique, Lyon 24-27 mai} Actes, Universit\'e Claude Bernard {\bf[1983]} 19. 

\vspace{4mm}\noindent[7] \textsc{J-F. Chamayou, J-L. Dunau}
{\it " Random difference equations with logarithmic distribution and the triggered shot noise "}
{\bf Adv. Appl. Math.} 29 {\bf[2002]} 454-470.
 
\vspace{4mm}\noindent[8] \textsc{J-F. Chamayou, J-L. Dunau}
{\it " Random difference equations: an asymptotical result "}
{\bf J. Comput. Appl. Math.} 154,1{\bf[2003]} 183-193.

\vspace{4mm}\noindent[9] \textsc{J-F. Chamayou}
{\it " Products of double gamma, gamma and beta distributions "}
{\bf Stat. Probab. Lett.} 68,2 {\bf[2004]} 199-208.


\vspace{4mm}\noindent[10] \textsc{E. A. Gussmann}
{\it Modifizierung der Gewichtsfunktionenmethode zur Berechnung der Fraunhoferlinien in Sonnen und Sternspektren }
{\bf Zeitschrift Astrophysik.} 65  {\bf[1967]} 456-497.

\vspace{4mm}\noindent[11] \textsc{J-L. Habemont, K. Hami-Eddine}
{\it " Simulation du bruit de grenaille  "}
{\bf Rapport 4GMM, Institut National Sciences Appliqu\'ees, mai} Universit\'e Toulouse{\bf[2005]}. 

\vspace{4mm}\noindent[12] \textsc{E. Hallen}
{\it "Iterated sine and cosine Integrals" }
{\bf Trans. Roy. Inst. Techn. Stockholm.} 12  {\bf[1947]} .

\vspace{4mm}\noindent[13] \textsc{E. R. Hansen } 
{\it A table of Series and Products }{\bf [1975]},Prentice Hall.

\vspace{4mm}\noindent[14] \textsc{K. S. Kolbig}
{\it On the integral $ \int_0^{\infty} e^{-\mu t} t^{\nu-1} log^{m} t dt $ }
{\bf Math. Comput.} 41 {\bf[1983]} 171-182.

\vspace{4mm}\noindent[15] \textsc{K. S. Kolbig}
{\it The Polygamma Function $ \psi^{(k)}(x) $ for $ x=\frac{1}{4} $ 
and $ x=\frac{3}{4} $  }
{\bf J. Comput. Appl. Math.} 75,1 {\bf[1996]} 43-46.

\vspace{4mm}\noindent[16] \textsc{V. Kourganoff, I.W. Busbridge}
{\it " Basic methods in Transfer Problems"}
{\bf[1952]} Clarendon, Oxford.

\vspace{4mm}\noindent[17] \textsc{H. Kunstner}
{\it "Zur Berechnung zweier spezieller Integrale" }
{\bf Wissenschõft. Zeit. Wilhem-Pieck Uni. Rostock, NaturW. Reihe.} 3  {\bf[1984]} 63-64.

\vspace{4mm}\noindent[18] \textsc{L. Lewin}
{\it "Polylogarithms and associated Functions"}
{\bf[1981]} North Holland, NewYork.

\vspace{4mm}\noindent[19] \textsc{G. Petiau}
{\it La Th\'eorie des Fonctions de Bessel }
{\bf CNRS.}   {\bf[1955]} .

\vspace{4mm}\noindent[20]\textsc{A.P. Prudnikov, Y.A. Brychkov, O.I.  Marichev} 
{\it Elementary Functions, Special Functions, More Special Functions, Vol.1 to  3 in Integrals and Series}{\bf [1992]}, Gordon Breach.

\vspace{4mm}\noindent[21] \textsc{V. G. Smith}
{\it An asymptotic expansion of $Ji_0(x)=\int_x^{\infty}\frac{J_0(t)}{t}dt$. }
{\bf J. Math. Phys. }22 {\bf [1943]} 58-59.

\vspace{4mm}\noindent[21] \textsc{J. Spanier, K.B.  Oldham}
{\it " An Atlas of Functions"}
{\bf[1987]} Hemisphere Publ, Washington D.C.; Springer Verlag, Berlin.

 \vspace{4mm}\noindent[22] \textsc{H. C. Van De Hulst}
{\it " Multiple Scattering"}
{\bf[1980]} Acad. Press, New York.

\vspace{4mm}\noindent[23] \textsc{K. Van Harn, F.W. Steutel}
{\it Infinite divisibility and the waiting-time paradox  }
{\bf Commun. Statist. Stochastic Models }11,3 {\bf [1995]} 527-540.

\newpage
\begin{center}
\begin{picture}(400,200)(0,0)

\put(200,0){\vector(1,0){100}}

\put(10,0){\line(0,1){200}}

\thicklines

\put(0,0){\line(1,0){200}}

\put(50,50){\oval(80,80)[bl]}

\put(100,180){\oval(100,100)[bl]}

\put(190,110){\oval(180,180)[bl]}
\put(260,190){\oval(140,140)[bl]}
\put(285,110){\oval(50,40)[bl]}

\thinlines
\put(145,190){$\Lambda_i$}
\put(1,10){...$t_{0}$}
\put(50,10){$t_{1}$}
\put(100,10){$t_{2}$...}
\put(160,10){...$t_{n-k}$...}
\put(260,10){...$t_{n}$}
\put(270,110){$W_t$}

\put(8,0){$\uparrow$}
\put(48,0){$\uparrow$}
\put(98,0){$\uparrow$}
\put(188,0){$\uparrow$}
\put(258,0){$\uparrow$}
\put(10,0){\line(0,1){200}}
\put(50,0){\line(0,1){200}}
\put(100,0){\line(0,1){200}}
\put(190,0){\line(0,1){200}}
\put(260,0){\line(0,1){200}}

\put(280,0){$t$}
\put(20,150){Scheme 1: Shot Noise Process}

\end{picture}

\end{center}

\newpage

\begin{center}
\begin{picture}(400,200)(0,0)

\put(200,0){\vector(1,0){100}}


\thicklines

\put(0,0){\line(1,0){200}}

\put(100,90){\oval(198,90)[bl]}

\put(190,150){\oval(180,180)[bl]}

\put(260,100){\oval(140,140)[bl]}

\put(281,200){\oval(40,40)[bl]}

\thinlines


\put(145,190){$\Lambda_i$}
\put(1,10){$t_{0}$}

\put(50,10){$t_{1}$}

\put(100,10){$t_{2}$}

\put(130,10){$t_{3}...$}
\put(160,10){...$t_{2(n-k)}$...}
\put(260,10){...$t_{2n}$}
\put(270,110){$W_t$}

\put(0,0){$\uparrow$}

\put(48,0){$\uparrow$}

\put(98,0){$\uparrow$}
\put(98,0){$\uparrow$}
\put(128,0){$\uparrow$}
\put(188,0){$\uparrow$}
\put(258,0){$\uparrow$}

\put(0,0){\line(0,1){200}}
\put(100,0){\line(0,1){200}}
\put(190,0){\line(0,1){200}}
\put(260,0){\line(0,1){200}}

\put(275,0){$t$}
\put(10,150){Scheme 2:Triggered Shot Noise Process}
\end{picture}
\end{center}

 \begin{figure}

\centerline{\includegraphics[width=15cm]{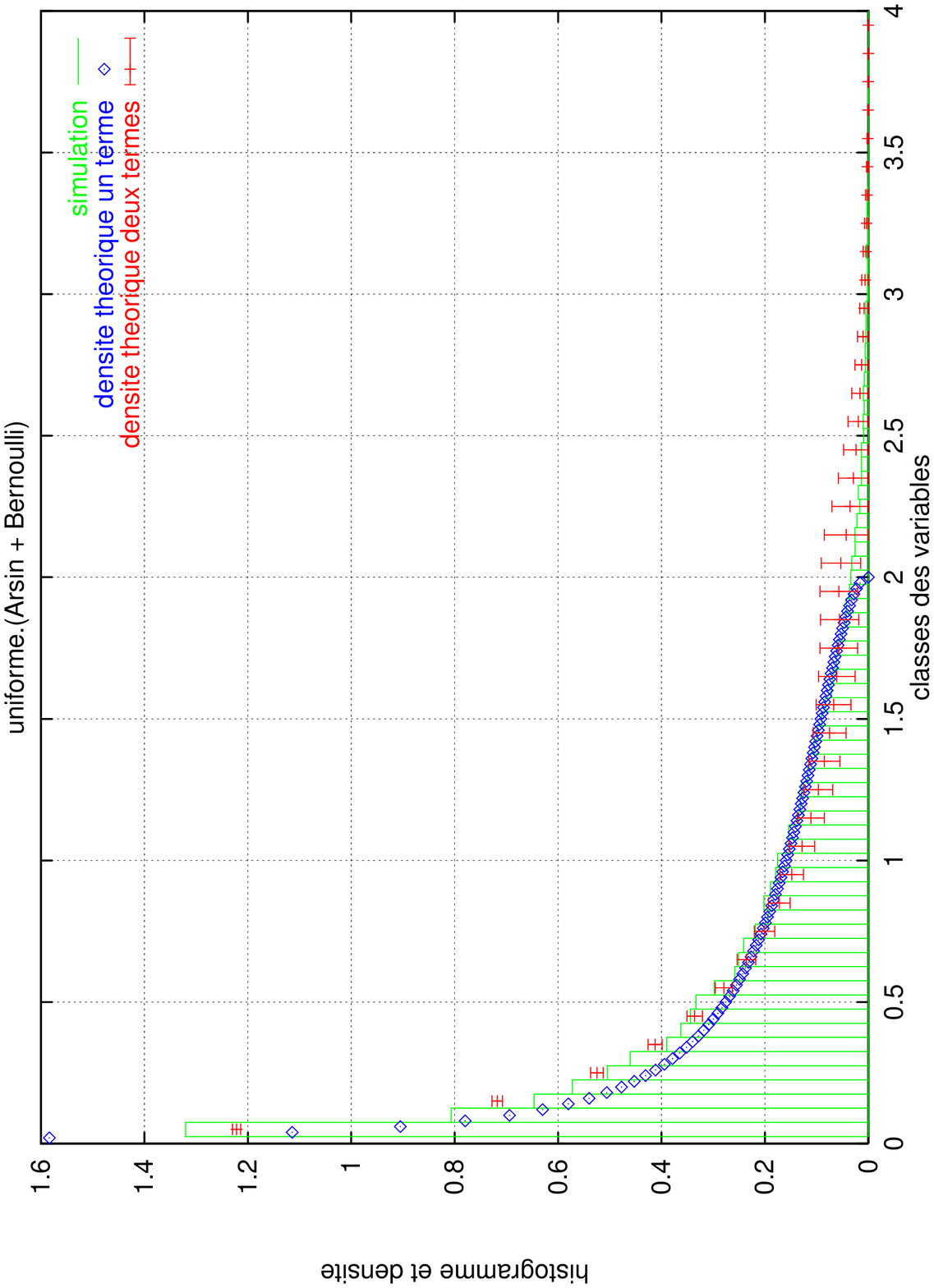}}
\caption{example}
\end{figure}
 \begin{figure}

\centerline{\includegraphics[width=15cm]{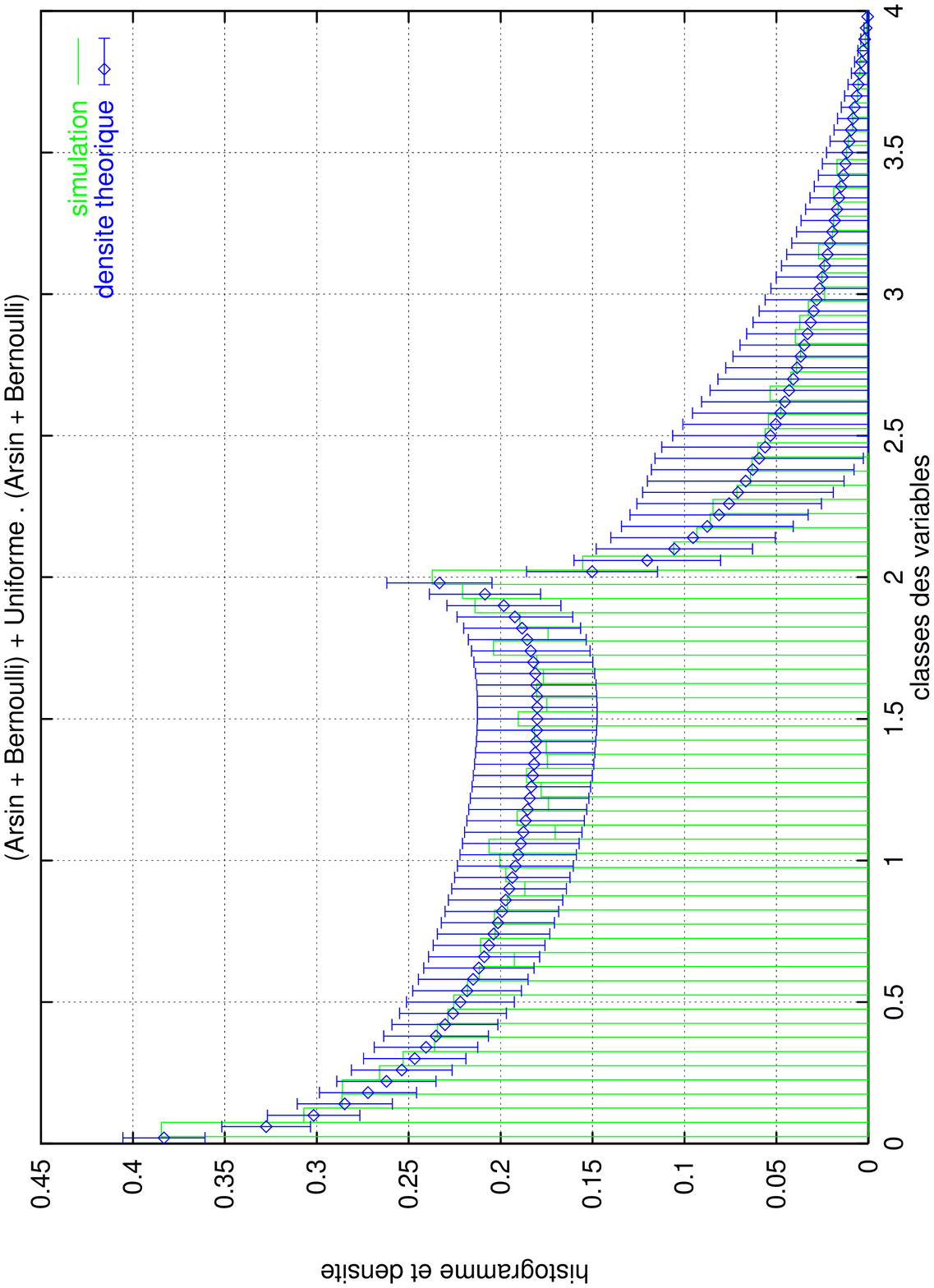}}
\caption{example}
\end{figure}
 \begin{figure}
\centerline{\includegraphics[width=15cm]{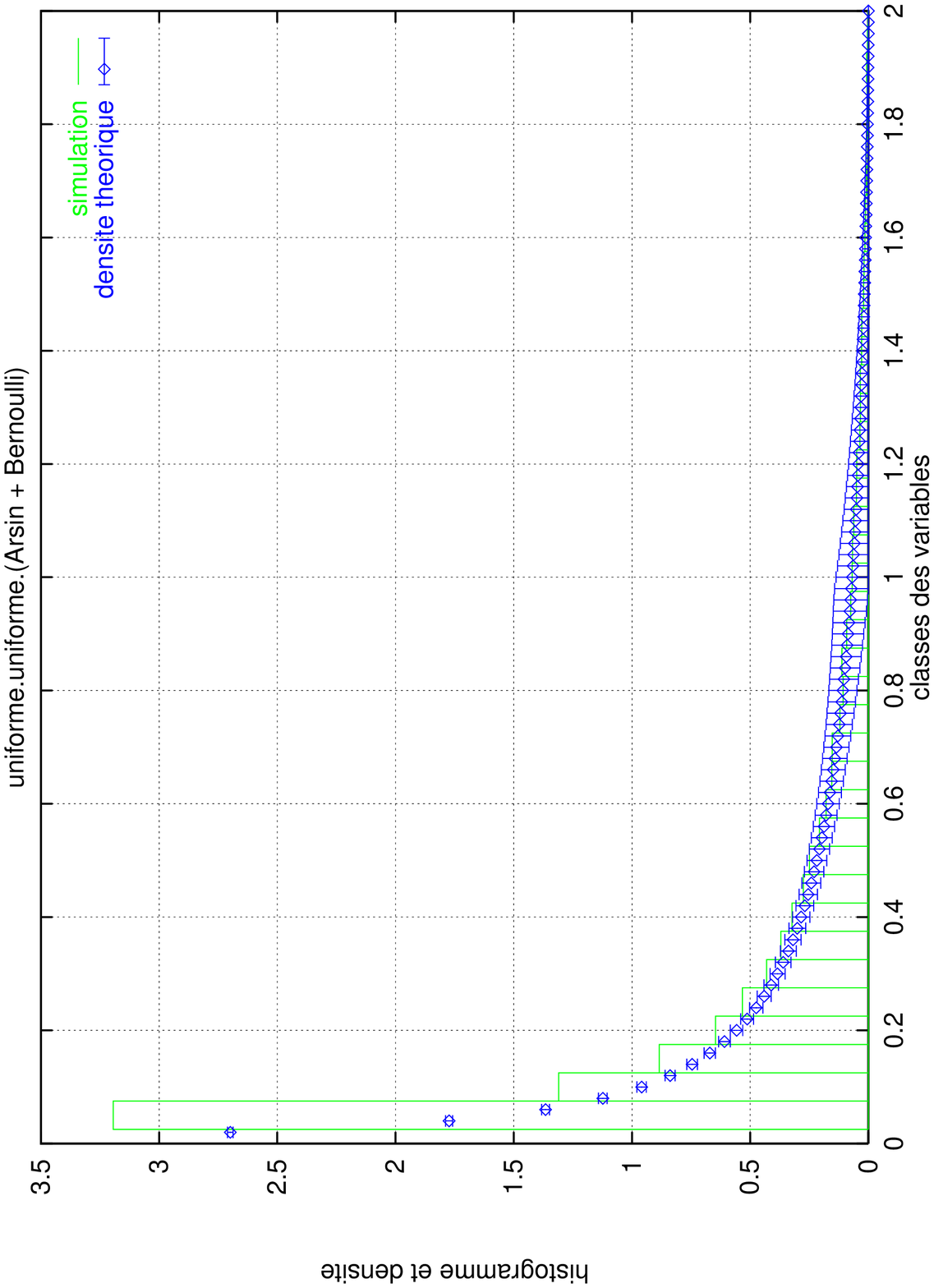}}

\caption{example}
\end{figure}
\begin{figure}
\centerline{\includegraphics[width=15cm]{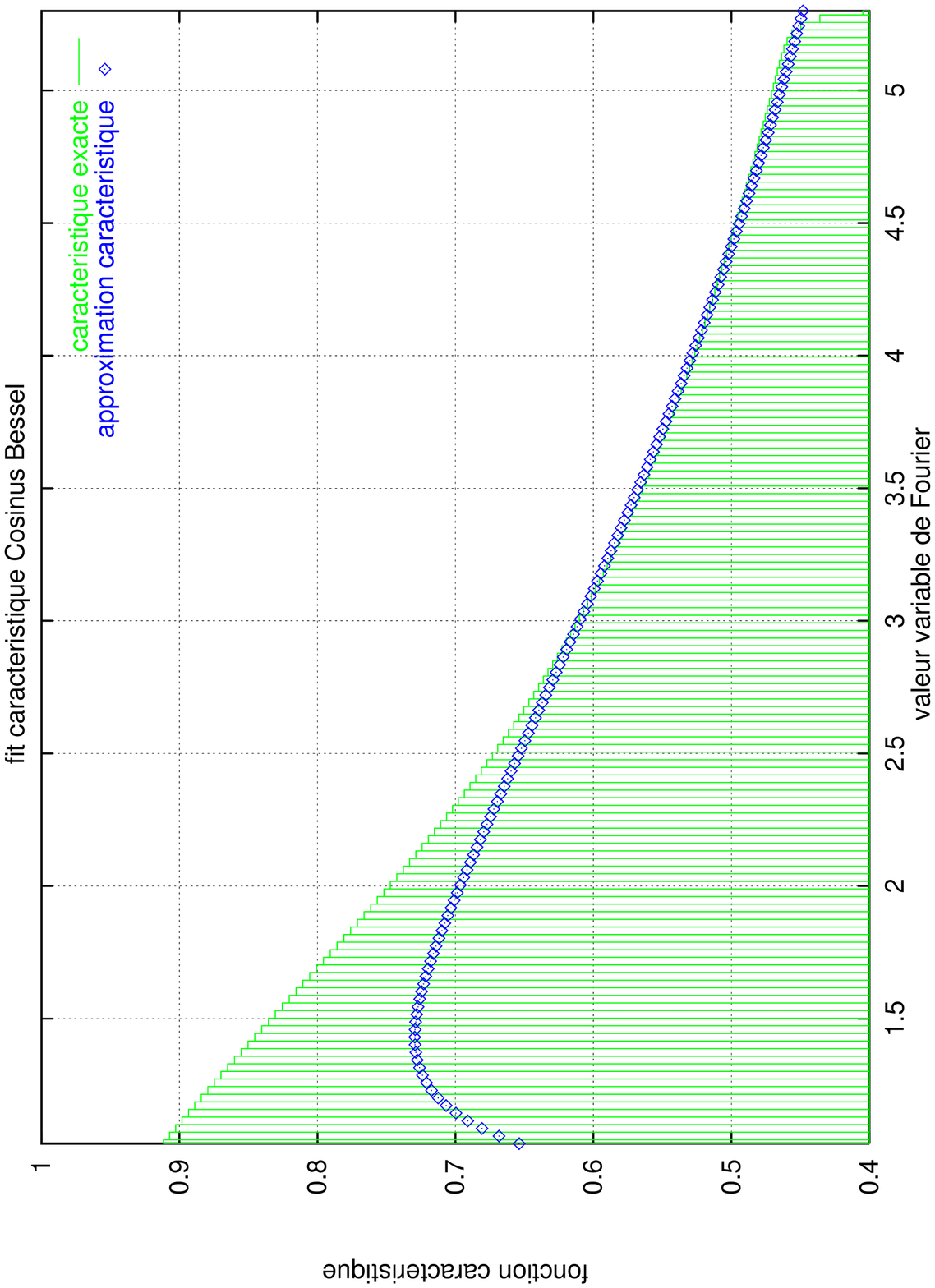}}
\caption{example}
\end{figure}

\begin{figure}
\centerline{\includegraphics[width=15cm]{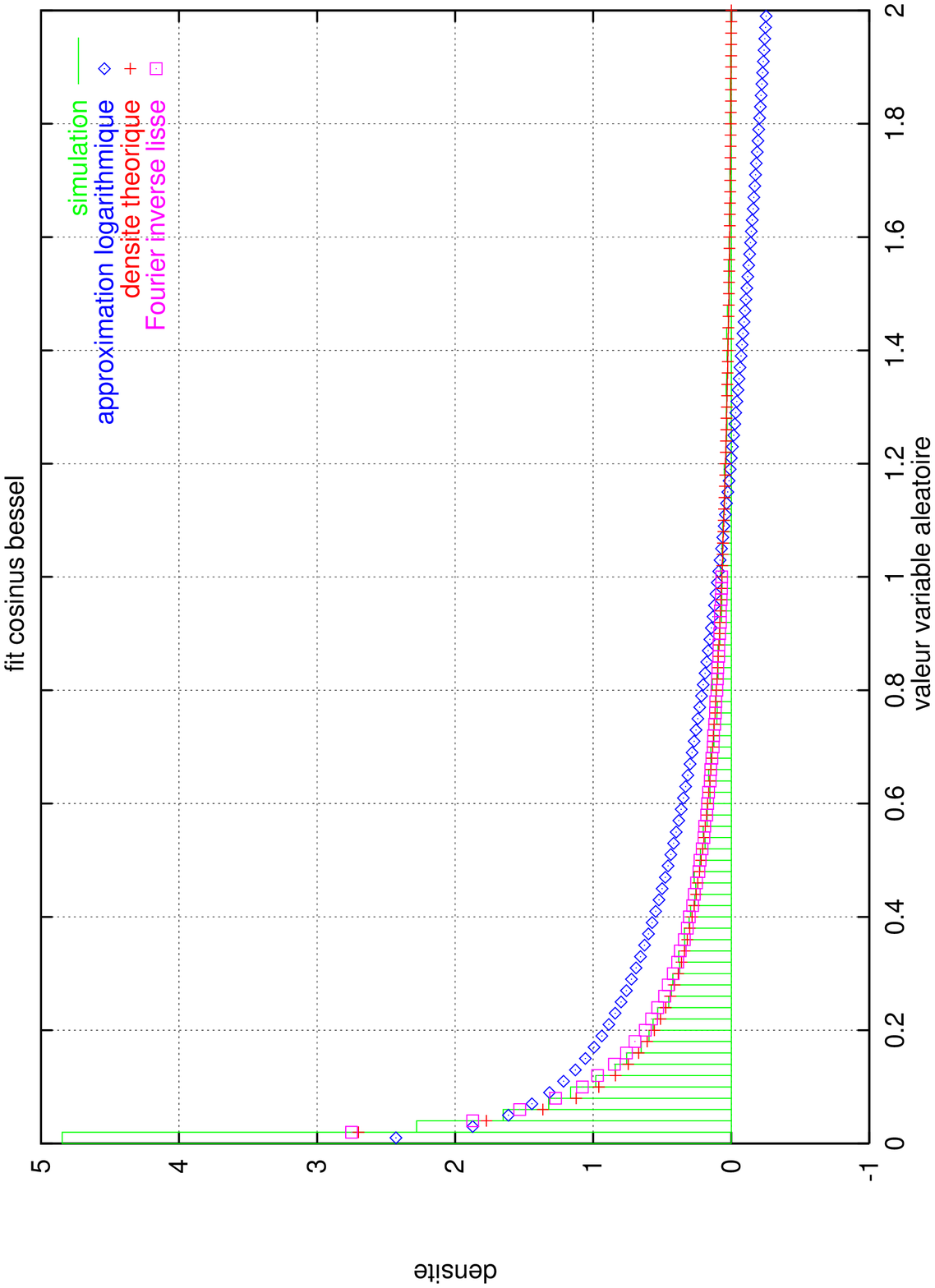}}
\caption{example}
\end{figure}

 \begin{figure}
\centerline{\includegraphics[width=15cm]{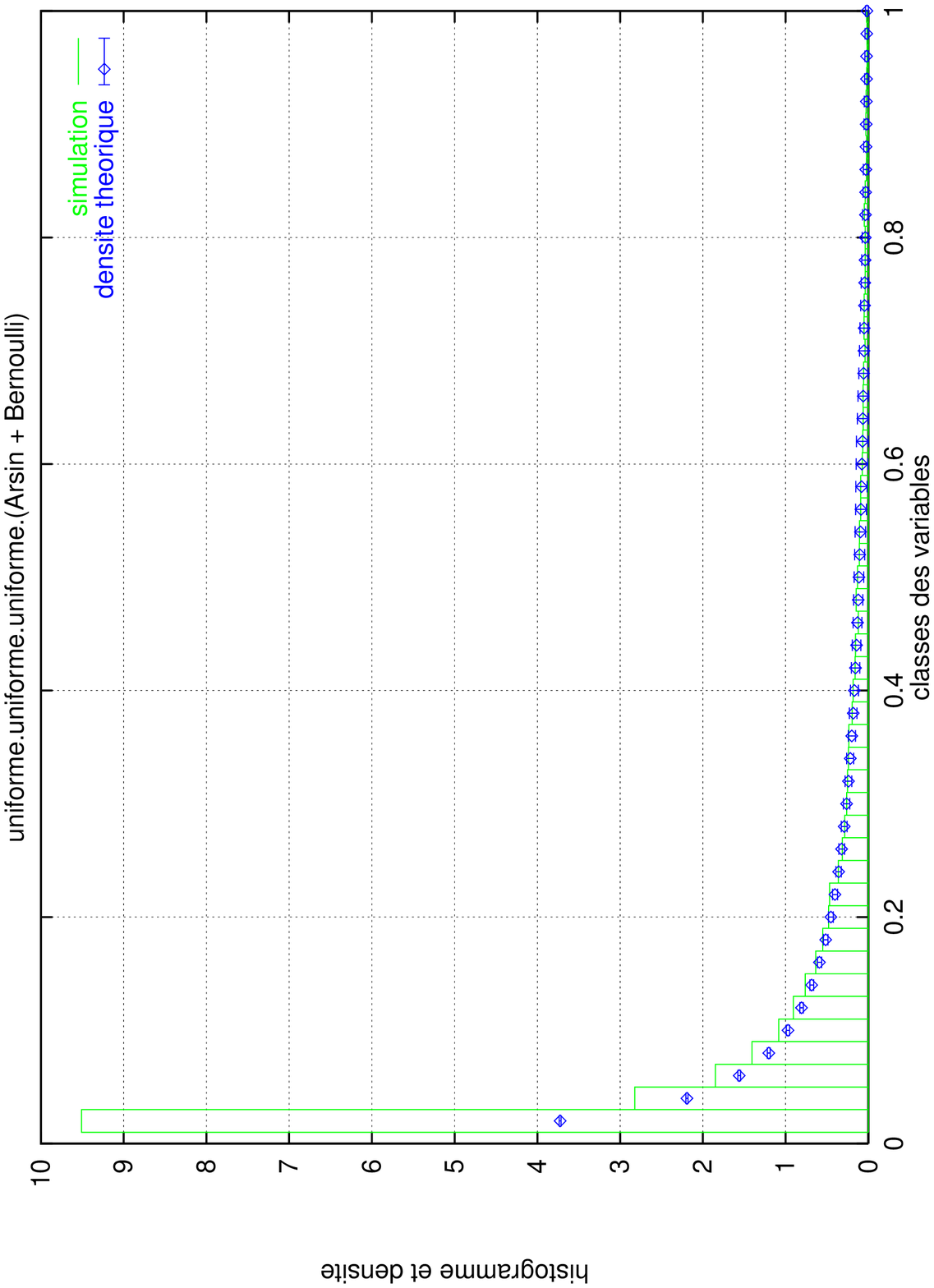}}
\caption{example}
\end{figure}

\end{document}